\documentclass[12pt]{amsart}
\usepackage[english]{babel}
\usepackage{amssymb,amsmath,amscd,amsthm}

\usepackage{cite}
\def\?{}

\newtheoremstyle{neosn}{0.5\topsep}{0.5\topsep}{\rm}{}{\sc}{.}{ }{\thmname{#1}\thmnumber{ #2}\thmnote{ {\mdseries#3}}}
\theoremstyle{neosn}

\newcommand{\Ann}{\,\mathrm{Ann}\,}

\newcommand{\Soc}{\,\mathrm{Soc}\,}

\newcommand{\ord}{\,\mathrm{ord}\,}

\begin{document}

\begin{center}
\large{\textbf{Rings which are Essential over their Centers, II}}
\end{center}

\hfill {\sf V.T. Markov}

\hfill Lomonosov Moscow State University

\hfill e-mail: vtmarkov@yandex.ru

\hfill {\sf A.A. Tuganbaev}

\hfill National Research University "MPEI"

\hfill Lomonosov Moscow State University

\hfill e-mail: tuganbaev@gmail.com

{\bf Abstract.} A ring $R$ with center
$C$ is said to be \textit{centrally
essential} if the module $R_C$ is an
essential extension of the module $C_C$.
We describe centrally essential exterior
algebras of finitely generated free
modules over not necessary commutative
rings and study properties of
semi-Artinian centrally essential rings.

V.T.Markov is supported by the Russian Foundation for Basic Research, project 17-01-00895-A.~ A.A. Tuganbaev is supported by Russian Scientific Foundation, project 16-11-10013.


\begin{center}
\textbf{1. Introduction}
\end{center}

All considered rings are associative unital rings.
The center of the ring $R$ is denoted by $C(R)$.

A ring $R$ with center $C$ is said to be \emph{centrally essential} if the module $R_C$ is an essential extension of the module $C_C$.

If $F$ is a field of characteristic $p\ge 0$ and $V=F^n$ is a finite-dimensional vector space over $F$ of dimension $n>0$,
then $\Lambda(V)$ denotes the exterior algebra of the space $V$ \cite[\S III.5]{Bourbaki}, i.e., $\Lambda(V)$ is a unital $F$-algebra with respect to the multiplication $\wedge$ with generators $e_1,\ldots,e_n$ and defining relations $e_i\wedge e_j=-e_j\wedge e_i$ for all $i,j\in\{1,\ldots,n\}$.

\textbf{Remark 1.1.} It is clear that any commutative ring and any direct factor of a centrally essential ring are centrally essential rings. There exist noncommutative finite centrally essential rings. For example, let $F$ be the field consisting of three elements, $V$ be a vector $F$-space with basis $e_1,e_2,e_3$, and let $\Lambda(V)$ be the exterior algebra of the space $V$. Since $e_1^2=e_2^2=e_3^2=0$ and any product of generators is equal to the $\pm$product of generators with ascending subscripts, $\Lambda(V)$ is a finite $F$-algebra of dimension 8 with basis
\begin{gather*}
\{1,e_1,e_2,e_3,e_1\wedge e_2,e_1\wedge e_3,e_2\wedge e_3,e_1\wedge e_2\wedge e_3\},\\
|\Lambda(V)|=3^8,\;e_k\wedge e_i\wedge e_j=-e_i\wedge e_k\wedge e_j=e_i\wedge
e_j\wedge e_k.
\end{gather*}
Therefore, if
$$
x=\alpha_0\cdot 1+\alpha_1^1e_1+\alpha_1^2e_2+\alpha_1^3e_3
+\alpha_2^1e_1\wedge e_2+\alpha_2^2e_1\wedge e_3+
\alpha_2^3e_2\wedge e_3+\alpha_3e_1\wedge e_2\wedge
e_3,$$ then
{\renewcommand{\arraycolsep}{0pt}
$$
\begin{array}{lll}[e_1,x]&=&2\alpha_1^2e_1\wedge
e_2+2\alpha_1^3e_1\wedge e_3,\\
{}[e_2,x]&=-&2\alpha_1^1e_1\wedge
e_2+2\alpha_1^3e_2\wedge e_3,\\
{}[e_3,x]&=-&2\alpha_1^2e_1\wedge
e_3-2\alpha_1^2e_2\wedge
e_3.\end{array}
$$
}
Thus, $x\in C(\Lambda(V))$ if and only if $\alpha_1^1=\alpha_1^2=\alpha_1^3=0.$ In other words, the center of the algebra $\Lambda(V)$ is of dimension 5. On the other hand, e.g., if $\alpha_1^1\neq 0$, then
$$
x\wedge(e_2\wedge e_3)=\alpha_0e_2\wedge e_3+\alpha_1^1e_1\wedge e_2\wedge e_3\in C(\Lambda(V))\setminus \{0\}.$$ In addition,
$e_2\wedge e_3\in C(\Lambda(V))$. Consequently, $\Lambda(V)$ is a
finite centrally essential noncommutative ring.

\textbf{Remark 1.2.} Let $F$ be a field of characteristic $p\ge 0$, $V$ be a finite-dimensional vector space over the field $F$, and let $\Lambda(V)$ be the exterior algebra of the space $V$. If $p=2$, then the ring $\Lambda(V)$ is commutative. In particular, $\Lambda(V)$ is a centrally essential ring. For $p\ne 2$, we can use the argument, which is similar to the above argument, to prove that $\Lambda(V)$ is a centrally essential ring
if and only if $\dim V$ is an odd integer.

In connection to Remark 1.2, we will prove Theorem 1.3 which is the first main result of the paper.

To formulate Theorem 1.3, we have to define the algebra $\Lambda(A^n)$ of the finitely generated free module $A^n$ of rank $n$ over a not necessarily commutative ring $A$ with center $C$. Namely,~ $\Lambda(A^n)=A\otimes_C\Lambda(C^n)$, where $\Lambda(C^n)$ is the exterior algebra of the free module $C^n$ over the commutative ring $C$; see \cite[\S III.5]{Bourbaki}.

Let $\{e_1,\ldots,e_n\}$ be a basis of the module $C^n$. By identifying $1\otimes x$ with $x$ for every $x\in \Lambda(C^n)$, we obtain that the set
$$
B_n=\{e_{i_1}\wedge\ldots\wedge e_{i_s}|0\leq s\leq n,\;\; 1\leq i_1<\ldots<i_s\leq n\}
$$
is a basis of the $A$-module $\Lambda(A^n)$ (we assume that the product is equal to 1 for $s=0$). It is directly verified that the ring $R=\Lambda(A^n)$ has a natural grading  $R=\oplus_{s\geq 0} R_s$, where $R_0=A$, $R_s=\oplus_{1\leq i_1<\ldots<i_s\leq n}Ae_{i_1}\wedge\ldots\wedge e_{i_s}$ for $1\leq s\leq n$, and $R_s=0$ for $s>n$.

\textbf{Theorem 1.3.} \textit{Let $A$ be a ring with center $C$ and $n$ a positive integer. The following conditions are equivalent}.\\
\textbf{1)} $R=\Lambda(A^n)$ is a \textit{centrally essential ring}.\\
\textbf{2)} $A$ is a \textit{centrally essential ring and at least one of the following conditions holds}:\\
$\phantom{xxx}$\textbf{a)} \textit{the ideal $\Ann_A(2)$ is an essential submodule of the module} $A_C$;\\
$\phantom{xxx}$\textbf{b)} $n$ \textit{is an odd integer.}

If $A$ is a ring of finite characteristic or $A$ does not have zero-divisors, then the formulation of Theorem 1.3 can be simplified; see Theorem 1.4 which is the main result of the paper.

\textbf{Theorem 1.4.} \textit{Let $A$ be a ring with center $C$ and
$n$ a positive integer.}

\textbf{1.} \textit{If $A$ is a ring of finite characteristic $s$ $($this is the case if the ring $A$ is finite$)$, then $\Lambda(A^n)$ is a centrally essential ring if and only if $A$ is a centrally essential ring and at least one of the following conditions holds:\\
$\phantom{xxx}$\textbf{a)} $s=2^m$ for some $m\in\mathbb{N}$;\\
$\phantom{xxx}$\textbf{b)} $n$ is an odd integer.}

\textbf{2.} \textit{If $A$ is a ring without zero-divisors, then
$\Lambda(A^n)$ is a centrally essential ring if and only if $A$ is a centrally essential ring and at least one of the following conditions holds:\\
$\phantom{xxx}$\textbf{a)} $A$ is a ring of characteristic $2$;\\
$\phantom{xxx}$\textbf{b)} $n$ is an odd integer.}

The third main result of the paper is Theorem 1.5.

\textbf{Theorem 1.5.} \textit{Let $R$ be a centrally essential ring with center $C$.}

\textbf{1.} \textit{If $R$ is a left or right semi-Artinian ring, then $R/J(R)$ is a commutative regular ring.}

\textbf{2.} \textit{If $R$ is a right nonsingular ring or $C$ is a semiprime ring, then $R$ is a commutative ring.}

\textbf{3.} \textit{If $B$ is a proper ideal of the ring $R$ generated by some central idempotents and the factor ring $A/B$ does not have nontrivial idempotents, then $A/B$ not necessarily is a centrally essential ring.}

\textbf{Remark 1.6.} In connection to Theorem 1.5(2), we remark that a ring with a semiprime center is not necessarily a semiprime ring. Indeed, the
ring $A$ of upper triangular $2\times 2$ matrices over a field $F$
is not a semiprime ring and the center $C$ of the ring $A$ is the
ring of scalar $2\times 2$ matrices over $F$ and hence $C$ is a
semiprime ring.

The proofs of Theorem 1.3, Theorem 1.4, and Theorem 1.5 are given
in the next section. We give some necessary notions.

Let $R$ be a ring with the center $C$. If $c\in C$, then we denote
by
$\Ann_R(c)$ the annihilator of the element $c$ in the ring $R$. We
denote by $J(R)$ the Jacobson radical of the ring $R$. For any two
elements $a,b\in R$, we set $[a,b]=ab-ba$. We denote by $\ord (a)$
the order of the element $a$ of the additive group $(R,+)$ of the
ring $R$.

A ring $R$ is said to be \emph{local} if $R/J(R)$ is a division ring.

A ring $R$ is said to be \emph{regular} if $a\in aRa$ for every element $a\in R$.

The \emph{socle} of the module $M$ is the sum $\Soc M$ of all its simple submodules; if $M$ does not have simple submodules, then $\Soc M=0$ by definition. A module $M$ is said to be \textit{semi-Artinian} if every its factor module is an essential extension of its socle.

A ring is said to be \textit{reduced} if it does not have non-zero nilpotent elements. A ring $R$ is said to be \textit{right nonsingular} if $R$ does not have non-zero elements whose right annihilators are essential right ideals.

All remaining necessary notions of ring theory can be found in the books \cite{Bourbaki,Lambek,Herstein}.

\begin{center}
\textbf{2. The proof of Theorems 1.3, 1.4, and 1.5}
\end{center}


\textbf{2.1. The proof of Theorem 1.3.} We set $R=\Lambda(A^n)$.

We assume that $R$ is a centrally essential ring.

Let $a\in A\setminus \{0\}$ and $a'=ae_1\wedge\ldots\wedge e_n$. Then $0\neq a'\in R$. Therefore, there exists an element $c\in C(R)$ such that $0\neq ca'\in C(R)$. We have $c=c_0+c'$, where $c_0\in A$, and $c'\in\oplus_{s>0}R_s$. It is easy to verify that $c_0\in C(A)=C(R)\cap R_0$. It is also clear that $ca'=c_0a'=c_0ae_1\wedge
\ldots \wedge e_n$, whence $c_0a\neq 0$. For every $b\in A$, we have
$$
0=[b,c_0ae_1\wedge\ldots\wedge e_n]=[b,c_0a]e_1\wedge\ldots
\wedge e_n,
$$
whence $c_0a\in C(A)$, i.e. $A$ is a centrally essential ring.

We assume that the ideal $\Ann_A(2)$ is not an essential submodule of the module $A_C$ and $n$ is an even integer.

We take an element $a\in A$ such that $a\neq 0$ and $Ca\cap\Ann_A(2)=0$. We consider the element $x=ae_2\wedge\ldots\wedge e_n$.

Let $c\in C(R)$ and $0\neq cx\in C(R)$. We have $c=c_0+c_1e_1+c'$,
where $c'$ is a linear combination of elements of the basis $B_n$ which are not equal to 1 and $e_1$. It is clear that $c_0,c_1\in C$ and $cx=c_0ae_2\wedge\ldots\wedge e_n+c_1ae_1\wedge\ldots\wedge e_n$, where the both summands are contained in the center of the ring $R$.

We prove that $c_0a=0$. Indeed, we have
\begin{multline*}
0=[e_1, c_0ae_2\wedge\ldots\wedge e_n]=c_0ae_1\wedge\ldots\wedge e_n-c_0ae_2\wedge\ldots\wedge e_n\wedge e_1=\\
c_0a(1-(-1)^{n-1})(-1)^{n-1}=2c_0ae_1\wedge\ldots\wedge
e_n,\end{multline*}
whence $c_0a\in \Ann_A(2)\cap Ca=0$ by the choice of $a$. Thus,
$c_1a\neq 0$ and $c_1e_1\in C(R)$. However, then $c_1\in C$ and $[c_1e_1,e_2]=2c_1e_1\wedge e_2=0$. Consequently, $c_1a\in Ca\cap\Ann_A(2)=0$. This is a contradiction.

Now we assume that $A$ is a centrally essential ring and at least one of the conditions \textbf{a)} or \textbf{b)} holds.

Let the condition \textbf{a)} hold. We set $N=\Ann_C(2)=C\cap\Ann_R(2)$. We remark that $N$ is an essential submodule in $A_C$. We consider an arbitrary non-zero element $x\in R$. We have
$$
x=\sum_{s=0}^n\sum_{1\leq i_1<\ldots<i_s\leq n}a_{i_1,\ldots,i_s}e_{i_1}\wedge\ldots\wedge e_{i_s},
$$
where the coefficients $a_{i_1,\ldots,i_s}$ are contained in $A$. We can multiply by elements $C\subseteq C(R)$ and obtain the situation where all the coefficients in the representation of $x$ are contained in $N$. Indeed, if some coefficient $a_{i_1,\ldots,i_s}$ is not contained in $N$, then there exists an element $c\in C$ such that $0\neq ca_{i_1,\ldots,i_s}\in N$, i.e. under the multiplication by $c$, the number of coefficients, not contained in $N$, decreases. It remains to remark that then $x\in C(R)$ provided all the coefficients in the representation of $x$ are contained in $N$. Indeed, $[x,a]=0$ for any $a\in A$, since $N\subseteq C(A)$ and
$$
[x,e_i]=\sum_{s=0}^n\sum_{i_1<\ldots<i_s}a_{i_1,\ldots,i_s}[e_{i_1}\wedge\ldots\wedge e_{i_s},e_i].
$$
We remark that if $s$ is even or $i\in \{i_1,\ldots,i_s\}$, then $[e_{i_1}\wedge\ldots\wedge e_{i_s},e_i]=0$. Otherwise, $[e_{i_1}\wedge\ldots\wedge e_{i_s},e_i]=\alpha e_{i_1}\wedge\ldots\wedge e_{i_s}\wedge e_i$, where $\alpha\in\{0,2\}$, i.e. we again have
$[a_{i_1,\ldots,i_s}e_{i_1}\wedge\ldots\wedge e_{i_s},e_i]=0$.
Since the elements of $A$ and $e_1,\ldots,e_n$ generate the ring $R$, we have $x\in C(R)$, which is required.

Now we assume that the condition \textbf{b)} holds. We consider an arbitrary non-zero element $x\in R$. By repeating the argument from the previous case, we can use the multiplication by elements of $C$ to obtain the situation such that all the coefficients $x$ with respect to the basis $B_n$ are contained in $C$.

We take the least odd $k$ such that the
element $e_{i_1}\wedge\ldots\wedge
e_{i_k}$ of the basis $B_n$ is contained
in the representation of $x$ with
non-zero coefficient $a$ (if this is not
possible, then $x\in C(R)$). Let $m=n-k$
and
$\{j_1,\ldots,j_q\}=\{1,\ldots,n\}\setminus
\{i_1,\ldots,i_k\}$. It is clear that the
integer $m$ is even, therefore,
$c=e_{j_1}\wedge\ldots\wedge e_{j_m}\in
C(R)$. Then it is easy to verify that
$cx=\pm ae_1\wedge\ldots\wedge e_n+x'$,
where \label{??}$x'$ is a linear
combination of elements of the basis
$B_n$ with even degree $s$ and
coefficients in $C$. Therefore, it
follows from the argument from the
previous case that $x'\in C(R)$. Finally,
it is directly verified that
$ae_1\wedge\ldots\wedge e_n\in C(R)$. The
assertion is proved.~\hfill$\square$

\textbf{Lemma 2.2.} \textit{If $A$ is a ring of finite characteristic $s$ and $C=C(A)$, then the following conditions are equivalent.}\\
\textbf{1)} \textit{The ideal $\Ann_A(2)$ is an essential submodule of the module} $A_C$.\\
\textbf{2)} $s=2^m$ for some $m\in\mathbb{N}$.

\textbf{Proof.} 1)\,$\Rightarrow$\,2).
We assume the contrary. Then there exists an odd prime integer $p$ dividing $s$. The non-zero ideal $\Ann_A(p)$ of the ring $A$ has the zero intersection with the ideal $\Ann_A(2)$. Therefore, the ideal $\Ann_A(2)$ is not an essential submodule of the module $A_C$. This is a contradiction.

2)\,$\Rightarrow$\,1). Since $s=2^m$, we have that for every $a\in A\setminus\{0\}$, the relation $\ord{a}=2^k$ holds for some $k\in \mathbb{N}$. Then $0\neq 2^{k-1}a\in Ca\cap\Ann_A(2)$. Therefore, the ideal $\Ann_A(2)$ is an essential submodule of the module $A_C$.~\hfill$\square$

\textbf{2.3. The proof of Theorem 1.4.} We set $R=\Lambda(A^n)$.

\textbf{1.} The assertion follows from Theorem 1.3 and Lemma 2.2.

\textbf{2.} If $A$ is a ring of characteristic $2$ or $n$ is an odd integer, then $R$ is a centrally essential ring by 1.

Now we assume that $A$ is a ring without zero-divisors and $R$ is a centrally essential ring. By Theorem 1.3, $A$ is a centrally essential ring and it is sufficient to consider the case, where $n$ is an even integer and the ideal $\Ann_A(2)$ is an essential submodule of the module $A_C$. Since $A$ is a ring without zero-divisors, $\Ann_A(2)=A$. Therefore, $A$ is a ring of characteristic $2$.~\hfill$\square$

\textbf{Proposition 2.4.}\label{idemp}
\textit{Let $R$ be a centrally essential ring with center $C$}.

\textbf{1.} \textit{In the ring $R$, all idempotents are central.}

\textbf{2.} \textit{If $R$ is a local ring, then the ring $R/J(R)$ is commutative.}

{\bf Proof.} \textbf{1.} Let $e\in R$ and $e^2=e$. It is sufficient to prove that $eR(1-e)=0=(1-e)Re$.

We assume that $a\in R$ and $ea(1-e)\ne 0$. Since $R$ is a centrally essential ring, $ea(1-e)c=d\ne 0$ for some central elements $c$ and
$d$. Then $0\ne ed(1-e)=de(1-e)=0$. This is a contradiction. Therefore, $eR(1-e)=0$. Similarly, $(1-e)Re=0$.

\textbf{2.} Let $x,y\in R$ and $xy-yx\not \in J(R)$. If $c\in C$ and $cx=d\in C\setminus \{0\}$, then $c(xy-yx)=dy-yd=0$, whence $c=0$. This is a contradiction.~\hfill$\square$

\textbf{Remark 2.5.} Let $R$ be a semiprime ring. It is well known that the set of all minimal right ideals of the ring $R$ coincides with the set of all minimal left ideals of the ring $R$ and this set coincides with the set of all right ideals $eR$ such that $eRe$ is a division ring. In addition, $\Soc R_R=\Soc _RR=S$.

\textbf{Remark 2.6.} Let $R$ be a ring and let $P$ be a nil-ideal of the ring $R$. It is well known that every idempotent $\overline e$ of the ring $R/P$ is of the form $e+P$, where $e=e^2\in R$.

\textbf{Proposition 2.7.} \textit{Let $R$ be a centrally essential ring, $P$ be a semiprime nil-ideal of the ring $R$, and let $\overline S$ be the right or left socle of the ring $R/J(R)$. If $\overline S$ is an essential left ideal of the ring $R/P$, then the ring $R/P$ is commutative.}

\textbf{Proof.} Let $h\colon R\to R/J(R)$ be the natural epimorphism. For every subset $X$ in $R$, we write $\overline X$ instead of $h(X)$. By Remark 2.5, there exists an ideal $S$ of the ring $R$ such that $P\subset S$ and $\overline S=\Soc _{\overline R}\overline R=\Soc \overline R_{\overline R}$.

First, we show that the ideal $\overline S$ is commutative. By Remarks 2.5 and 2.6, any minimal left ideal $V$ of the ring $\overline R$ is generated by some primitive idempotent $\overline e$ which is of the form $\overline e=e+J(R)$ for some primitive idempotent $e$ of the ring $R$. By Proposition 2.4, all idempotents of the ring $R$ are central. If $0\neq x\in \overline R\overline e$, then it follows from the minimality of the left ideal $V$ that
there exists an element $y\in \overline R$ such that $yx=\overline e= (y\overline e)x$, i.e., $V$ is a division ring. Since $J(Re)=Re\cap J(R)$, we have that $Re$ is a local ring and $R=Re\oplus R(1-e)$. Therefore, $R$ is a centrally essential ring and the division ring $V=Re/J(Re)$ is commutative by Proposition 2.4.
Since $S=\Soc(\overline R)$ is a direct sum of some set of minimal ideals, each of them is a field, we have that the ideal $\overline S$ is commutative.

Since $\overline S$ is an essential left ideal of the semiprime ring $\overline R$, the right annihilator and the left annihilator of the ideal $S$ are equal to the zero. For any two elements $x,y\in \overline R$ and arbitrary $a,b\in \overline S$, we have
\begin{multline*}
ab[x,y]=b(xy-yx)a=(bx)(ya)-(by)(xa)=\\
=yabx-byxa=byax-abyx=abyx-abyx=\overline 0.
\end{multline*}
Then $\overline{S}^2[x,y]=0$, whence $[x,y]=\overline 0$. Therefore, $R/P$ is a commutative ring.~\hfill$\square$

\textbf{Proposition 2.8.} \textit{Let $A$ be a centrally essential ring with center $C$. The following conditions are equivalent.}

\textbf{1)} $A$ is a \textit{semiprime ring};

\textbf{2)} $C$ is a \textit{semiprime ring};

\textbf{3)} $A$ is a \textit{reduced ring};

\textbf{4)} $A$ is a right nonsingular ring;

\textbf{5)} $A$ is a \textit{commutative reduced ring}.

\textbf{Proof.} The implication 5)\,$\Rightarrow$\,1) is ovbious.

1)\,$\Rightarrow$\,2). Let $c\in C$ and  $c^2=0$. Then $(AcA)^2=0$. Since $A$ is a semiprime ring, $C$ is a reduced ring. In particular, $C$ is a semiprime ring.

2)\,$\Rightarrow$\,3). We assume the contrary. Let $a$ be a
non-zero element of the ring $A$ such that $a^2=0$. Since $A$ is a
centrally essential ring, there are two non-zero elements $c,d\in
C$ such that $ac=d\ne 0$. Then $d^2=(ac)^2=a^2c^2=0$, whence
$(CdC)^2=0$. This is a contradiction.

3)\,$\Rightarrow$\,4). The assertion is true for any reduced ring.

4)\,$\Rightarrow$\,2). Let $x\in A$ and $x^2=0$. Let us assume
that $x\ne 0$. There exist central elements $c,d\in C$ with
$xc=d\ne 0$. By the Zorn lemma, there exists a right ideal $Y$ of
the ring
$A$ such that $xcA\cap Y=0$ and  $xcA\oplus Y$ is an essential
right ideal. Then $d(xA)=xcxA=x^2cA=0$. In addition, $dY=Yd\in
xcA\cap Y=0$. Therefore, $d(xcA\oplus Y)=0$. In addition,
$xcA\oplus Y$ is an essential right ideal. Since the ring $A$ is
right nonsingular, $d=0$. This is a contradiction.

3)\,$\Rightarrow$\,5). We verify that for any $r\in A$, the ideal $$r^{-1}C=\{c\in C:\;rc\in C\}$$ is dense in $C$. Indeed, let $d\in C$ and $dr^{-1}C=0$. If  $dr=0$, then $d\in r^{-1}C$ and $d^2=0$, whence
$d=0$. Otherwise, since $A$~is a centrally essential ring, there exists an element $z\in C$ such that $zdr\in C\setminus\{0\}$. Then $zd\in r^{-1}C$  and $(zd)^2=0$, whence $zd=0$; this contradicts to the choice of $z$. The property `$r(r^{-1}C)\neq 0$ for every $r\in A\setminus\{0\}$' is equivalent to the property that $A$ is a centrally essential ring.
Therefore, $A$~is a right ring of quotients of the ring $C$ in the sense of \cite[\S 4.3]{Lambek}. Then $A$ can be embedded in the complete ring of quotients of the ring $C$ which is commutative.~\hfill$\square$

\textbf{2.9. The proof of Theorem 1.5.}

\textbf{1.} Let $R$ be a centrally essential right or left semi-Artinian ring and $\overline R=R/J(R)$. Since $R$ is a right or left semi-Artinian ring, $J(R)$ is a nil-ideal \cite[Proposition 3.2]{NP}. By Proposition 2.7, $R/J(R)$ is a commutative ring. Every commutative semi-Artinian semiprimitive ring is regular by \cite[ Theorem 3.1]{NP}.

\textbf{2.} The assertion follows from Proposition 2.8.

\textbf{3.} Let $F$ be the field consisting of three elements, $A=\Lambda(F^3)$ be the exterior algebra of the three-dimensional vector $F$-space $F^3$, and let $S=\Lambda(F^2)$ be the exterior algebra of the two-dimensional vector $F$-space $F^2$ which is considered as a subalgebra of the algebra $A$. We consider the direct product $P=A^{\mathbb N}=\{(a_1,a_2,\ldots)\,|\,a_i\in A\}$
of the countable set of copies of the ring $A$ and its
subring $R$ consisting of all finally constant sequences
$(a_1,a_2,\ldots)\in P$ which are stabilized at elements of $S$ at
finite step depending on the sequence. Let $e_i$ be a central idempotent which has the identity element of the field $F$ on the  $i$th position and zeros on the remaining positions. We denote by $B$ the ideal of the ring $R$ generated by all the idempotents
$\{e_i\}$. It follows from Theorem 1.4 that $R$ is a centrally essential ring and the factor ring $R/B$ is isomorphic to the ring
$S$ which is not centrally essential and does not have nontrivial idempotents.~\hfill$\square$

In connection to Remark 1.1, Theorem 1.5(1),(2) and Proposition
2.4(2), we formulate the following open questions.

\textbf{2.10. Open questions.} Let $R$ be a centrally essential ring with center $C$, Jacobson radical $J(R)$ and prime radical $N(R)$.

\textbf{1.} Is it true that the ring $R/J(R)$ is commutative?

\textbf{2.} Is it true that the ring $R/N(R)$ is commutative?

\textbf{3.} If the ring $R$ is semiperfect, is it true that $\Soc(R_C)=\Soc(R_R)$?

\textbf{4.} If the ring $R$ is semiperfect, is it true that $R=C+J(R)$?

\end{document}